\documentclass[10pt]{article}

\usepackage{amsmath}
\usepackage{amssymb}
\usepackage{graphicx}
\usepackage{epic,eepic}

\usepackage{color}
\renewcommand{\div}{\text{div}}

\usepackage{amstext,amsthm,amssymb,amsmath}
\usepackage{graphicx}
\usepackage{subfigure}
\usepackage{wrapfig}
\usepackage{fullpage}
\usepackage{color}
\usepackage{multirow}
\usepackage{tabulary}
\usepackage{booktabs}
\usepackage{enumerate}

\newtheorem{theorem}{Theorem}[section]

\newtheorem{remark}[theorem]{Remark}

\newcommand{\pa}{\partial}

\newcommand{\norm}[1]{\left\|#1\right\|}

\newcommand{\forma}{a(\, \cdot \, , \, \cdot \,)}
\newcommand{\dualforma}{a^\ast(\, \cdot \, , \, \cdot \,)}

\hfuzz=\maxdimen
\title{Goal-oriented adaptivity for GMsFEM}

\author{
Eric T. Chung\thanks{Department of Mathematics, The Chinese University of Hong Kong, Hong Kong SAR. Eric Chung's research
is supported by the Hong Kong RGC General Research Fund project 400813.},\and
  Wing Tat Leung\thanks{Department of Mathematics, Texas A\&M University, College Station, TX },\and 
 Sara Pollock\thanks{Department of Mathematics, Texas A\&M University, College Station, TX }
}

\begin{document}
\maketitle

\begin{abstract}

In this paper we develop two goal-oriented adaptive strategies for 
{\em a posteriori} error estimation within the generalized multiscale 
finite element framework. 
In this methodology, one seeks to determine the number of multiscale basis 
functions adaptively for each coarse region 
to efficiently reduce the error in the goal functional. 
Our first error estimator uses a residual based strategy
where local indicators on each coarse neighborhood are the product of local
indicators for the primal and dual problems, respectively.
In the second approach, viewed as the multiscale extension of the 
dual weighted residual method (DWR), the error indicators are computed 
as the pairing of the local $H^{-1}$ residual of the primal problem
weighed by a projection into the primal space of the $H_0^1$ dual solution 
from an enriched space, over each coarse neighborhood.
In both of these strategies, the goal-oriented indicators 
are then used in place of a standard
residual-based indicator to mark coarse neighborhoods of the mesh for
further enrichment in the form of additional multiscale basis functions.
The method is demonstrated on high-contrast problems with 
heterogeneous multiscale coefficients, and is seen to outperform 
the standard residual based strategy with respect to efficient reduction
of error in the goal function.

\end{abstract}

\section{Introduction}
\label{sec:intro}

Many practical problems are multiscale in nature,
including flow in porous media, seismic wave
propagation, and physical processes in perforated media.
These problems are described
by partial differential equations (PDEs) with potentially high contrast
multiscale coefficients. 
Direct computation of high resolution discrete solutions to these problems 
can be very expensive. 
Typically, some type of model reduction techniques are used to solve 
multiscale problems. 
Established techniques include numerical homogenization methods
\cite{dur91, cdgw03, cd07} 
and multiscale methods \cite{hw97, ehw99, egh12, chung2014adaptive, 
chan2015adaptive, hkj12, jennylt03, ce09,chung2015generalizedperforated, 
ArPeWY07,AKL, Efen_GVass_11}. 
In numerical homogenization methods, the upscaled media
properties are computed over coarse grid blocks,  
each of which is much larger than a characteristic length scale.
In multiscale methods, local multiscale basis functions determined by
local fine-scale problems are constructed in 
each element of the coarse grid. Generally, one uses a few multiscale basis 
functions in each coarse element to 
approximate the global solution by solving a coarse mesh problem 
over the entire domain.
Multiscale basis functions are constructed in an offline step before the 
coarse mesh problem is solved, after which some type of adaptivity
is needed to choose multiscale basis function appropriately.

In this paper we will use multiscale methods, where multiscale basis functions 
are constructed in each coarse region, as
illustrated schematically in Figure~\ref{schematic}. 
To be more specific, we consider a multiscale problem
\[
L(u) = f,
\]
where $L=-\div(\kappa(x)\nabla u)$ and seek the solution in the form
\[
u(x)=\sum_{i,j} u_{i,j} \phi_{i}^{\omega_j}.
\]
In each coarse region $\omega_j$, we construct a set of multiscale basis functions $\phi_i^{\omega_j}$,
$i=1,...,N_j$. These multiscale basis functions, as described below, will be constructed in the offline
stage using the generalized multiscale finite element method (GMsFEM)
 and represent the local heterogeneities of the solution space. 
In earlier works on the multiscale finite element method (MsFEM) \cite{hw97},
the authors sought one multiscale basis function per coarse element. 
However as it was later argued,
one may need additional basis functions over each coarse
element for a sufficiently high fidelity approximation.
It is further shown that the optimal number of multiscale basis functions 
in each region depends on the heterogeneities in the solution space. 
Typically, some adaptive criteria based on {\em a posteriori} error estimation
is used to determine how many basis 
functions to choose in each coarse region $\omega_i$. 
For exampe, in \cite{chung2014adaptive},
the authors develop an error indicator based on the $H^{-1}$ norm
of the residual to determine the number of basis 
functions to add in each region over each adaptive iteration.
Adaptive multiscale methods follow traditional adaptivity concepts as in 
\cite{verfurth96,BDD04,CKNS08,NSV09,Hols2001a,MeNo05,Stevenson07}; 
however, the multiscale indicators contain information about local 
heterogeneities. Earlier approaches to goal-oriented adaptive methods
for multiscale problems
include numerical regularization and numerical homogenization with adaptive
mesh refinement as in \cite{OdVe00a,JhDe12a}.  To the authors' knowledge,
 the current presentation is the first to develop a goal-oriented enrichment 
strategy within the general GMsFEM framework. 

For many practical problems, one is interested in approximating some
function of the solution, known as the quantity of interest, rather than
the solution itself. 
Examples include an average or weighted
average of the solution over a particular 
subdomain, or some localized solution response. In these cases,
goal-oriented adaptive methods yield a more efficient approximation
than standard adaptivity, as the enrichment of degrees of freedom is 
focused on the local improvement of the quantity of interest rather than
across the entire solution 
\cite{PrOd99,BeRa01,MoSt09,EHL02,Giles.M;Suli.E2003,Gratsch.T;Bathe.K2005,BaRa03,BuNa15,HPZ11a}. 
In this paper, we study goal-oriented adaptivity for multiscale methods, and
in particular the design of error indicators to drive the adaptive enrichment
based on the goal function.
In multiscale methods, goal-oriented adaptivity can play an important role 
in the efficient approximation of the quantity of interest as
heterogeneities in the coefficients may require standard adaptive methods to 
add degrees of freedom in regions with limited influence on the goal function.
In this paper, we develop a goal-oriented approach for multiscale methods 
within GMsFEM framework.
In the proposed approach, we increase the accuracy of the approximation
by enriching the space rather than refining the mesh
by choosing multiscale basis functions computed in the offline stage.

For multiscale basis construction, we use GMsFEM. The construction of 
multiscale basis functions uses local snapshot spaces and requires solving 
local spectral problems over each coarse element. 
The local snapshot functions represents
the solution space in each coarse region, and they can include all possible 
local 
fine-grid functions or harmonic functions. In the snapshot space, we perform 
a local spectral decomposition and select multiscale basis functions 
which correspond to the dominant eigenvalues. The multiscale basis functions 
are constructed by multiplying the dominant eigenmodes by a partition of
unity function, {\em e.g.,} multiscale partition of unity function. 
In \cite{chung2014adaptive}, we developed an adaptive approach and developed 
{\em a posteriori} error indicators, which include the information from local 
spectral problems, {\em e.g.,} the value of the eigenvalue corresponding to the 
first eigenvector not included in the coarse space. 
We derived error estimates and presented numerical results
which demonstrate the improved efficiency of the adaptive approach, 
guided these indicators.
For goal-oriented problems, we now  design goal-oriented 
error indicators, which are different from those developed earlier 
\cite{chung2014adaptive} for multiscale problems, by the additional 
consideration of a dual problem to direct the adaptivity towards the
approximation of the quantity of interest.

In this paper we develop two goal-oriented adaptive strategies for 
{\em a posteriori} error estimation. 
Our first error estimator uses an idea similar to a standard residual based 
adaptive method, and can be seen as the multiscale extension of the 
$hp$-adaptive method presented in~\cite{BuNa15}.
In this case the elementwise indicator is formed by the product of local
residual indicators for the primal and dual problems, respectively.
In the second approach, viewed as the multiscale extenstion of the 
dual weighted residual method (DWR), the error indicators are computed 
as the pairing of the local $H^{-1}$ residual of the primal problem
weighed by a projection into the primal space of the $H_0^1$ dual solution 
from an enriched space.
In both of these strategies, the goal-oriented indicators are then 
used in place of a standard
residual-based indicator to mark coarse elements of the mesh for
further enrichment in the form of additional multiscale basis functions.

The remainder of the paper is organized as follows.  In Section \ref{prelim}
we present and overview of multiscale methods and adaptivity.  In Section
\ref{cgdgmsfem} we give a detailed description of the 
construction of the multiscale basis functions, review some results on 
residual-based adaptivity for GMSFEM.  In Section 
\ref{sec:errorindicator} we introduce two goal-oriented
{\em a posteriori} error indicators and present an algorithm for
goal-oriented adaptivity. Finally in Section \ref{sec:numresults} we present
numerical results demonstrating the efficiency of the proposed 
goal-oriented error indicators.

\section{Overview of Concepts}
\label{prelim}
In this paper, we consider second order multiscale elliptic problems of the form
\begin{equation} \label{eq:original}
\begin{split}
-\mbox{div} \big( \kappa(x) \, \nabla u  \big) &=f \quad \text{in} \quad D, \\
u &= 0 \quad \text{on} \quad \partial D,
\end{split}
\end{equation}
where $D$ is the computational domain, and 
$\kappa(x)$ is a scalar valued heterogeneous coefficient
with multiple scales and high contrast.
The problem (\ref{eq:original}) can be solved by many classical numerical 
techniques, such as the conforming finite element method,
but with extremely high computational complexity
due to the fact that a very fine mesh is necessary to resolve the multiscale 
nature of the solution. 
Thus, some multiscale model reductions are needed to compute an accurate 
solution efficiently.
In the following, we give a brief overview of GMsFEM and its basis 
enrichment techniques as applied to problem \eqref{eq:original}.
Let $u\in V = H^1_0(\Omega)$ be the true solution satisfying
\begin{equation}\label{eqn:primal}
a(u,v) = (f,v), \quad v\in V,
\end{equation}
where
$ \displaystyle a(u, v) = \int_D \kappa(x) \nabla u \cdot \nabla v \, dx$, and $ \displaystyle (f,v) = \int_D f v \, dx$.
Define the energy norm on $V$ by  $\|u\|_V^2 = a(u,u)$.

\begin{figure}[htb]
  \centering
  \includegraphics[width=0.65 \textwidth]{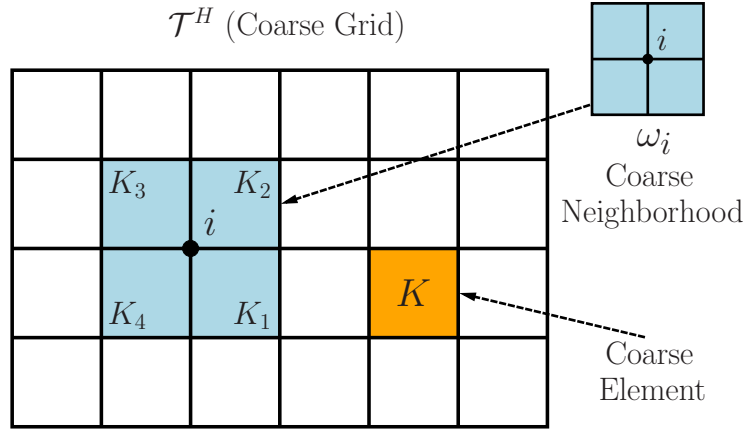}
  \caption{Illustration of a coarse neighborhood and a coarse element}
  \label{schematic}
\end{figure}

To introduce the GMsFEM for the problem (\ref{eq:original}), we first give 
the notion of fine and coarse grids.
We let $\mathcal{T}^H$ be a standard conforming triangulation of the computational domain $D$ into finite elements,
which can be triangular, rectangular or some other polygons. 
We refer to this partition as the coarse grid.
Subordinate to  the coarse grid, 
we define the fine grid partition, denoted by $\mathcal{T}^h$,
by refining each coarse element into a connected union of fine grid blocks. 
We assume the above refinement is performed such that 
$\mathcal{T}^h$ is a conforming partition of $D$. 
We let $N$ be the number of interior coarse grid nodes, and 
let $\{x_i\}_{i=1}^{N}$ be the set of coarse grid nodes or vertices of
the coarse mesh $\mathcal{T}^H$. Moreover, we define the coarse neighborhood of the node $x_i$ by
\begin{equation} \label{neighborhood}
\omega_i=\bigcup\{ K_j\in\mathcal{T}^H; ~~~ x_i\in \overline{K}_j\},
\end{equation}
which is the union of all coarse elements which have the node $x_i$ as a vertex.
See Figure~\ref{schematic} for an illustration of the coarse elements and  
coarse neighborhoods within the coarse grid.
We emphasize the use of $\omega_i$ to denote a coarse neighborhood, 
and $K$ to denote a coarse element throughout the paper.

Next, we briefly overview the continuous Galerkin (CG) formulation of GMsFEM,
a generalization of the classical MsFEM \cite{hw97}.
For each coarse node $x_i$, we define a set of basis functions supported 
on the neighborhood $\omega_i$.
We denote the $k$-th basis function supported on the coarse neighborhood
$\omega_i$ by $\psi_k^{\omega_i}$, 
We remark that in the GMsFEM, 
we will use multiple basis functions per coarse neighborhood,
and the index $k$ represents the local numbering of these basis functions.
These multiscale basis functions are constructed from a local snapshot space
and a local spectral decomposition defined on that snapshot space. 
The snapshot space contains a collection of many basis functions that can 
be used to capture most of the fine features of the solution,
and the multiscale basis functions $\psi_k^{\omega_i}$ are constructed 
by selecting the dominant modes of a local spectral problem. 
Using the these multiscale basis functions, the CG solution is represented 
as $u_{\text{ms}}(x)=\sum_{i,k} c_{k}^i \psi_{k}^{\omega_i}(x)$.
Once the basis functions are identified, the CG global coupling is given through the variational form
\begin{equation}
\label{eq:globalG} a(u_{\text{ms}},v)=(f,v), \quad \text{for all} \, \, v\in
V_{\text{off}},
\end{equation}
where  $V_{\text{off}}$ is the space spanned by the basis functions $\{ \psi_k^{\omega_i} \}$,
and $a(\cdot,\cdot)$ is the usual bilinear form corresponding to (\ref{eq:original}).
We remark that one can use other formulations, such as the 
discontinuous Galerkin formulation (see e.g., \cite{Wave,WaveGMsFEM,eglmsMSDG}),
the mixed formulation (see e.g., \cite{chung2015mixed,chan2015adaptive})
or the hybridized discontinuous Galerkin formulation (see e.g., \cite{efendiev2015multiscale})
to couple the multiscale basis functions.

In using GMsFEM, it is desirable to determine
the number of basis functions per coarse neighhorbood
adaptively based on the heterogeneities of the coefficient $\kappa(x)$
in order to obtain an efficient representation of the solution. 
In \cite{chung2014adaptive}, a residual based {\em a posteriori} error indicator 
is derived and an adaptive basis enrichment algorithm is developed under the CG formulation. 
In particular, it is shown that 
\begin{equation*}
\| u - u_{\text{ms}}\|_V^2 \leq C \sum_{i=1}^N r_i^2,
\end{equation*}
where $r_i$ is the residual of the solution $u_{\text{ms}}$ on the coarse neighborhood $\omega_i$.
Thus, local residuals of the multiscale solution
give indicators to the error of the solution in the energy norm,
and one can add basis functions to the coarse neighborhoods
when the residuals $r_i$ are large. 
Convergence of this adaptive basis enrichment algorithm is also shown in \cite{chung2014adaptive}.
On the other hand, for some applications one needs to adaptively construct 
new basis functions in the online stage
in order to capture distant effects.
In \cite{chung2015residual}, such online adaptivity is proposed and mathematically analyzed.
More precisely, when the local residual $r_i$ is large, one can construct 
a new basis function $\phi \in V_0(\omega_i)$ in the online stage by solving
\begin{equation*}
a(\phi, v) = (r_i,v), \quad \forall v \in V_0(\omega_i),
\end{equation*}
where $V_0(\omega_i)$ is the restriction of $V$ in $\omega_i$
with zero trace on $\partial\omega_i$.
Numerical results in \cite{chung2015residual} show that a couple of these online basis functions
can help to reduce the error by a large amount.

The adaptivity procedures discussed above are designed with the aim of 
reducing the error in the energy norm. 
In some applications, one may be more interested in reducing error measured
by some function of the solution other than a norm.
For example, in flow applications, one needs to obtain a good approximation
of the pressure in locations where the wells are situated. 
Therefore, we now consider goal-oriented adaptivity within GMsFEM. 
Specifically, we define a linear functional
$g:V \rightarrow \mathbb{R}$, referred to as the goal functional. 
In goal-oriented adaptivity, one wants to adaptively enrich the approximation space
in order to reduce the goal error defined by $g(u-u_{\text{ms}})$. 
In the construction of goal-oriented adaptivity for GMsFEM, 
we use local indicators based on the solution of a dual problem:
finding $z\in V$ such that
\begin{equation}\label{eqn:dual}
a^*(z,v) = g(v), \quad \forall v\in V.
\end{equation}
For a primal problem $a(u, v) = (f,v)$ based on bilinear form
$\forma$, the dual form $\dualforma$ is the formal adjoint of the primal,
satisfying $a^\ast(w,v) = a(v,w)$, and in 
the current symmetric case, $\dualforma$ is identical to the primal.
Formally, the primal-dual equivalence follows for $u$ the solution to the
primal problem~\eqref{eqn:primal} and $z$ the solution to the dual
problem~\eqref{eqn:dual}
\begin{align}\label{eqn:pd_equiv}
f(z) = a(u,z) = a^\ast(z,u) = a(z,u) = g(u).
\end{align}
Error estimates for the quantity of interest $g(u)$ follow 
from~\eqref{eqn:pd_equiv} and Galerkin orthogonality with respect to the
discrete problems and their respective solutions.
Forming error indicators based on both primal and dual problems and these
estimates, 
we add multiscale basis functions to coarse 
neighborhoods when the values of the corresponding indicators are large. 
Our numerical examples show that the goal-oriented approach performs better than
the residual approach for high-contrast problems
when the error is measured by the goal-functional, $g(u-u_{\text{ms}})$.

\section{The GMsFEM and residual-based adaptivity}
\label{cgdgmsfem}

In this section, we will give a detailed description
of the GMsFEM (see for example \cite{egh12, eglp13}) and it's residual based adaptivity 
(see for example \cite{chung2014adaptive}). 

\subsection{Local basis functions}
\label{locbasis}
We first present the construction
of the multiscale basis functions.
This construction is performed in the offline stage;
that is, basis functions are pre-computed before the actual solve
of the problem. 
The construction starts with a snapshot space. 
This space contains a relatively large set of basis functions which
can be used to capture most features of the fine-scale solution. 
The next step is to perform a local dimension reduction to obtain a 
lower dimensional subspace that can still be used to approximate the 
solution with good accuracy. 
The local dimension reduction is performed by solving a spectral problem,
and the dominant eigenfunctions are used as the  multiscale basis functions.

First, we define a snapshot space $V_{\text{snap}}^{\omega_i}$,
where the functions in $V_{\text{snap}}^{\omega_i}$ are supported in $\omega_i$. 
The snapshot space can be the space of all fine-scale basis functions $V(\omega_i) = \{ v |_{\omega_i} \, | \, v\in V\}$
or the solutions of some local problems with various choices of boundary conditions.
For example, we can use the following $\kappa$-harmonic extensions to form a 
snapshot space. 
Specifically, let $\{x_j^i\}, ~j = 1, \ldots L_i$, 
index the set of fine-grid vertices that lie on the boundary of each 
coarse neighborhood, $\pa \omega_i$. 
Define the unit source 
functions $\delta_j^h(x) = \delta(x_j^i)$ for each $j = 1, \ldots L_i$.
Then construct the snapshot function $\psi_{j}^{\omega_i, \text{snap}} \in V(\omega_i)$ by solving
\begin{equation} \label{harmonic_ex}
\begin{split}
-\text{div}(\kappa(x) \nabla \psi_{j}^{\omega_i, \text{snap}} ) &= 0,
\quad \text{in} \, \, \, \omega_i, \\
\psi_{j}^{\omega_i, \text{snap}} &= \delta_j^h, \quad \text{on} \, \, \, \partial\omega_i.
\end{split}
\end{equation}
The snapshot space $V_{\text{snap}}^{\omega_i}$ corresponding to the 
region $\omega_i$, then contains $L_i$ functions 
$$
V_{\text{snap}}^{\omega_i} = \text{span}\{ \psi_{j}^{ \omega_i,\text{snap}}:~~~ 1\leq j \leq L_i \}.
$$
We define the corresponding change of variable matrix 
$$
R_{\text{snap}}^{i} = \left[ \psi_{1}^{\omega_i,\text{snap}}, \ldots, \psi_{L_i}^{\omega_i,\text{snap}} \right],
$$
where $\psi_{j}^{\omega_i,\text{snap}}$ are considered as the columns of the matrix. 

We next determine a set of dominant modes from $V_{\text{snap}}^{\omega_i}$,
and the resulting lower dimensional space is called the offline space $V_{\text{off}}^{\omega_i}$.
To construct the offline space $V_{\text{off}}^{\omega_i}$, we perform a 
dimension reduction of the space of snapshots using an auxiliary spectral 
decomposition. The analysis in \cite{egw10} motivates the following generalized 
eigenvalue problem for eigenvalues $\lambda_k^{\text{off}}$ and
eigenfunctions $\Psi_k^{\text{off}}$ in the space of snapshots:
\begin{equation} \label{offeig}
A^{\text{off}} \Psi_k^{\text{off}} = \lambda_k^{\text{off}} S^{\text{off}} \Psi_k^{\text{off}},
\end{equation}
where
\begin{equation*}
 \displaystyle A^{\text{off}} = [a_{mn}^{\text{off}}] = \int_{\omega} \kappa(x) \nabla \psi_m^{\text{snap}} \cdot \nabla \psi_n^{\text{snap}} = (R^{i}_{\text{snap}})^T A R^i_{\text{snap}},
 \end{equation*}
 \begin{center}
 and
 \end{center}
 \begin{equation*}
 \displaystyle S^{\text{off}} = [s_{mn}^{\text{off}}] = \int_\omega  \widetilde{\kappa}(x)\psi_m^{\text{snap}} \psi_n^{\text{snap}} = (R^i_{\text{snap}})^T S R^i_{\text{snap}},
\end{equation*}
where $A$ and $S$ denote analogous fine-scale stiffness and mass
matrices as defined by
\begin{equation*}
A_{ij} = \int_{D} \kappa(x) \nabla \phi_i \cdot \nabla \phi_j \, dx,
\quad
S_{ij} = \int_{D} \widetilde{\kappa}(x)  \phi_i  \phi_j \, dx,
\end{equation*}
where $\phi_i$ is the fine-scale basis function for $V$. We will give the 
definition of $\widetilde{\kappa}(x)$ later on.
To generate the offline space we then select the smallest $l_i$ eigenvalues 
from Equation~\eqref{offeig} and form the corresponding eigenvectors in the 
space of snapshots by setting
$\psi_k^{\omega_i,\text{off}} 
= \sum_{j=1}^{L_i} \Psi_{kj}^{\text{off}} \psi_j^{\omega_i,\text{snap}}$ 
(for $k=1,\ldots, l_i$), where $\Psi_{kj}^{\text{off}}$ are the 
coordinates of the vector $\Psi_{k}^{\text{off}}$, and $l_i$ is the number
of eigenvectors chosen to span the offline space.
We will use the set $\{ \psi_k^{\omega_i,\text{off}} \}$ 
of local basis functions to form the approximation space in the next section.
\subsection{CG formulation}
\label{globcoupling}

In this section 
we create an appropriate solution space and the variational formulation for a 
continuous Galerkin approximation of Equation~\eqref{eq:original}. 
The idea is to use the basis set 
$\{ \psi_k^{\omega_i,\text{off}} \}_{k=1}^{l_i}$, $i=1,2,\cdots, N$,
to form the approximation space, called the offline space, and apply the 
standard continuous Galerkin formulation.
We begin with an initial coarse space 
$V^{\text{init}}_0 = \text{span}\{ \chi_i \}_{i=1}^{N}$, where we recall 
$N$ denotes the number of coarse neighborhoods corresponding to interior 
coarse nodes. 
Here, $\chi_i$ are the standard multiscale partition of unity functions 
which are supported in $\omega_i$ and are defined by
\begin{align} \label{pou}
-\text{div} \left( \kappa(x) \, \nabla \chi_i  \right) &= 0, 
\quad \text{in} \, \, \, K \subset \omega_i, \\ 
\chi_i &= g_i, \quad \text{on} \, \, \, \partial K \backslash\partial\omega_i, 
\nonumber  \\
\chi_i &= 0, \quad \text{on} \, \, \, \partial\omega_i, \nonumber
\end{align}
for all coarse elements $K \subset \omega_i$, where $g_i$ is a continuous 
function on $\partial K$ which is linear on each edge of $\partial K$.
Based on the analysis in \cite{egw10}, the summed, pointwise energy 
$\widetilde{\kappa}$ required for the 
eigenvalue problems (\ref{offeig}) is defined as
\begin{equation*}
\widetilde{\kappa} = \kappa \sum_{i=1}^{N} H^2 | \nabla \chi_i |^2,
\end{equation*}
where $H$ denotes the coarse mesh size.

The partition of unity functions $\chi_i$ are then multiplied by the 
eigenfunctions $\{ \psi_k^{\omega_i,\text{off}} \}_{k=1}^{l_i}$
to construct the multiscale basis functions
\begin{equation} \label{cgbasis}
\psi_{i,k} = \chi_i \psi_k^{\omega_i, \text{off}}, \quad \text{for} \, \, \,
1 \leq i \leq N \, \, \,  \text{and} \, \, \, 1 \leq k \leq l_i,
\end{equation}
where we recall $l_i$ denotes the number of offline eigenvectors 
that are chosen for each coarse node $i$. 
We note the construction in Equation~\eqref{cgbasis} yields  continuous 
basis functions due to the multiplication of offline eigenvectors with the 
initial (continuous) partition of unity. 
Next, we define the continuous Galerkin spectral multiscale space as
\begin{equation} \label{cgspace}
V_{\text{off}}  = \text{span} \{ \psi_{i,k} : \,  \, 1 \leq i \leq N \, \, \,  \text{and} \, \, \, 1 \leq k \leq l_i  \}.
\end{equation}
%
Using this offline space, we obtain the GMsFEM
as in (\ref{eq:globalG}).



\subsection{Residual-based adaptivity}

In this section, we give a brief review of the adaptive basis enrichment 
algorithm proposed in \cite{chung2014adaptive}. After solving the 
coarse mesh problem and computing error indicators, the standard
D\"orfler marking strategy is applied with respect to the neighborhoods 
indexed by vertices $i = 1, \ldots, N$, rather than the elements $K$.  
The marked neighborhoods are then enriched with additional basis functions. 
The spectral problem (\ref{offeig}) gives a natural ordering
of basis functions in the local snapshot space $V^{\omega_i}_{\text{snap}}$
with respect to the eigenvalues, in increasing order.
The analysis of \cite{chung2014adaptive} then suggests an adaptive procedure to 
add basis functions based on a local error indicator. 
The error indicator is defined by a $H^{-1}$-norm based residual,
which gives a robust indicator with good performance for cases with 
high contrast media. 

Let $\omega_i$ be a coarse neighborhood and write $V_i = V_0(\omega_i)$.
For a given multiscale solution $u_{\text{ms}}$, 
we define a linear functional $R^u_i(v)$ on $V_i$ by
\begin{equation}\label{eqn:primalresi}
R^u_i(v) =  \int_{\omega_i} fv - \int_{\omega_i} \kappa \nabla u_{\text{ms}}\cdot \nabla v, \quad v\in V_i.
\end{equation}
The norm of $R^u_i$ is defined as
\begin{equation}
\| R^u_i \|_{V_i^*} = \sup_{v\in V_i} \frac{ |R^u_i(v)| }{\| v\|_{V_i}}.
\label{eq:Rnorm}
\end{equation}
where $\|v\|_{V_i}^2 =\int_{\omega_i} \kappa(x) | \nabla v|^2 \, dx$.
The norm $\| R^u_i \|_{V_i^*}$ gives a measure on how well the solution $u_{\text{ms}}$ 
satisfies the variational problem (\ref{eqn:primal}) restricted to $V_i$.
The norm $\| R^u_i \|_{V_i^*}$ can be obtained by solving an auxiliary problem of finding $w\in V_i$ such that 
$a(w,v) = R_i^u(v)$ for all $v\in V_i$ and then setting $\| R^u_i \|_{V_i^*}  = \|w\|_V$.
One can solve this auxiliary problem in the snapshot space $V_{\text{snap}}^{\omega_i}$ to reduce the computational cost.
Moreover, in \cite{chung2014adaptive} 
the following {\em a posteriori} error bound is proved
\begin{equation}
\label{eq:a-post-bound}
\| u-u_{\text{ms}}\|_V^2 \leq C_{\text{\rm err}}  \sum_{i=1}^N \|R^u_i\|^2_{V_i^*}  (\lambda^{\omega_i}_{l_i+1})^{-1},
\end{equation}
where $C_{\text{err}}$ is a uniform constant independent of the contrast of $\kappa$,
and $\lambda_j^{\omega_i}$ is the $j$-th eigenvalue
of the eigenvalue problem (\ref{offeig}) for the coarse neighborhood $\omega_i$.
In particular, $\lambda_{l_i+1}^{\omega_i}$ is the first, {\em i.e.,}
smallest eigenvalue from the spectral problem \eqref{offeig} for which the
corresponding eigenvector is not included in the construction of the offline 
space.

Using the above error bound (\ref{eq:a-post-bound}), a convergent adaptive 
enrichment algorithm is developed in \cite{chung2014adaptive}.
We now describe this algorithm. 
The algorithm is an iterative process, and basis functions are added in 
each iteration/level based on the magnitudes of local residuals. 
We use $m \geq 1$ to index the enrichment level
and let $V_{\text{off}}^m$ be the solution space at the $m$-th iteration.
For each coarse region, let $l_i^m$
be the number of eigenfunctions used at the enrichment level $m$
for the coarse region $\omega_i$.

\subsection{Adaptive algorithm}\label{subsec:adaptive_algorithm}
The adaptive enrichment algorithm for GMsFEM is now summarized below.
Choose a fixed marking parameter $0 < \theta < 1$.
Choose also an initial offline space $V_{\text{off}}^1$ by specifying a fixed 
number of basis functions for each coarse neighborhood,
and this number is denoted by $l_i^1$, for each $i = 1, \ldots, N$.
Then, generate a sequence of spaces $V_{\text{off}}^m$
and a sequence of multiscale solutions $u_{\text{ms}}^m$ obtained by 
solving \eqref{eq:globalG}.
Specifically, for each $m=1,2,\cdots$,  perform the following calculations:

\begin{enumerate}
\item[Step 1:] Find the multiscale solution in the current space. That is,
find $u_{\text{ms}}^m \in V^m_{\text{off}}$ such that
\begin{equation}
a(u^m_{\text{ms}}, v) = (f, v) \quad \text{for all} \,\,\, v \in V^m_{\text{off}}.
\label{eq:solve}
\end{equation}

\item[Step 2:] Compute the local residual. For each coarse region $\omega_i$,
compute
\[
\eta^2_i = \|R_i^u\|^2_{V_i^*}  (\lambda^{\omega_i}_{l^m_i+1})^{-1},
\]
where
\begin{equation*}
R_i^u(v) =  \int_{\omega_i} fv - \int_{\omega_i} \kappa(x) 
\nabla u^m_{\text{ms}}\cdot \nabla v,
\end{equation*}
consistent with \eqref{eqn:primalresi},
and the norm is defined in (\ref{eq:Rnorm}) respectively.
Next, re-enumerate the coarse neighborhoods so the above local residuals 
$\eta_i^2$ are arranged in decreasing order
$\eta^2_1 \geq \eta^2_2 \geq \cdots \geq \eta^2_N$.
That is, in the new enumeration,
the coarse neighborhood $\omega_1$ has the largest residual $\eta_1^2$
and the coarse neighborhood $\omega_N$ has the least residual $\eta_N^2$.
\begin{remark}\label{remark:bin}
An alternate approach to avoid the $N \log N$ complexity of the full sort
is the standard binning or heapifying strategy \cite{MoSt09}. 
Let $\eta^2 = \sum_{i = 1}^n \eta_i^2$, and consider only $\eta_i$ that
satisfy $\eta_i^2 > (1 - \theta) \eta^2/N$.  Let $M = \max_i \eta_i^2$,
and perform a partial sort of the remaining indicators collecting or binning
the indices for which $2^{-p}M \le \eta_i^2 < 2^{-(p+1)}M$, for 
$p = 0, 2, \ldots q$, where $q$ is the smallest integer to satisfy
$2^{-(q+1)}M \le (1 - \theta) \eta^2/N.$
\end{remark}

\item[Step 3:] Find the coarse regions where enrichment is needed. 
Choose the smallest integer $k$ such that
\begin{equation}
\theta \sum_{i=1}^N \eta_i^2 \leq \sum_{i=1}^k \eta_i^2.
\label{eq:criteria1}
\end{equation}
The coarse neighborhoods $\omega_1,\omega_2,\cdots, \omega_k,$ are then
enriched with additional basis functions. If the partial sort of 
Remark~\ref{remark:bin} is used in place of the sort, elements are marked by 
emptying the first bin, those indicators with $M \le \eta_i^2 < M/2$, and then
continuing on to the second bin, and so forth until~\eqref{eq:criteria1} is 
satisfied.  As elements within bins are not sorted, this yields a quasi-optimal
marked set, {\em i.e.,} the marked set may not be the set of least-cardinality
to satisfy~\eqref{eq:criteria1} as in the full sort, but it is within a 
factor of two of the least cardinality.

\item[Step 4:] Enrich the space. For each $i=1,2,\cdots, k$, add basis function
for the region $\omega_i$ according to the following rule.
Let $s$ be the positive integer such that
$\lambda_{l_i^m+s+1}$ is large enough compared with $\lambda_{l_i^m+1}$
(see Remark \ref{remark:large}).
Then include the eigenfunctions 
$\Psi^{\text{off}}_{l_i^m+1}, \cdots, \Psi^{\text{off}}_{l_i^m+s}$
in the construction of the basis functions.
The resulting space is denoted as $V_{\text{off}}^{m+1}$.
Mathematically, the space $V_{\text{off}}^{m+1}$ is defined as
\begin{equation*}
V_{\text{off}}^{m+1} = V_{\text{off}}^{m} + \text{span} \cup_{i=1}^k \cup_{j=l_i^m+1}^{l_i^m+s} \{  \psi_{i,j} \}
\end{equation*}
where $\psi_{i,j} = \chi_i \psi_j^{\omega_i,\text{off}}$
and $\psi_j^{\omega_i,\text{off}} = \sum_{r=1}^{l_i} \Psi_{jr}^{\text{off}} \psi_r^{\text{snap}}$,
with $j=l_i^m+1, \cdots, l_i^m+s$,
denote the new basis functions obtained by the eigenfunctions $\Psi^{\text{off}}_{l_i^m+1}, \cdots, \Psi^{\text{off}}_{l_i^m+s}$.
In addition, we set $l_i^{m+1} = l_i^m+s$.

\begin{remark}
\label{remark:large}
The mathematical analysis in \cite{chung2014adaptive} specifies 
the choice of $s$. 
In practice, one can take $s=1$ since the eigenvalues in (\ref{offeig}) 
have fast growth. 
\end{remark}

\end{enumerate}

\section{Goal-oriented adaptivity}
\label{sec:errorindicator}


In this section, we present a goal-oriented adaptive enrichment 
algorithm for GMsFEM. 
The goal-oriented variant of the adaptive method requires the solution of 
a dual problem in addition to the primal at each iteration. The indicators are
computed with both the primal residual and either the dual residual or 
a projection of an enriched dual solution into the primal space.  These 
indicators predict which neighborhoods to enrich to increase the 
quality of the approximation of the quantity of interest.
After introduction of the discrete dual problem, the finite dimensional
analogue of~\eqref{eqn:dual}, 
we propose two error indicators for goal-oriented enrichment.

The dual problem plays a vital role in goal-oriented adaptivity as 
the vehicle for introducing the goal functional $g$ into the adaptive process. 
Given a goal functional 
$g:V \rightarrow \mathbb{R}$, we define the discrete dual problem 
on approximation space $V_{\text{off}} \subset V$ as: 
find $z \in V_{\text{off}}$ such that
\begin{equation}
a(v,z) = g(v), \;\forall v\in V_{\text{off}}.
\label{eqn:discrete_dual}
\end{equation}

As in~\eqref{eqn:dual},
the discrete dual form $\dualforma$ is identical to the primal $\forma$
for symmetric  problems. 
The dual problem, however, features the goal functional $g$ as the source.
The discrete dual solution, may now be used 
to define goal-oriented error indicators.

\subsubsection*{$H^{-1}$-based goal-oriented indicator}
Our first goal-oriented indicator is similar in form to the  
residual based indicator described in the previous section. 
To motivate this indicator, we introduce the local bilinear form
$a(u,v)_i = \int_{\omega_i}\kappa(x)  \nabla u \cdot \nabla v \, dx$, 
and the induced localized energy norm $\norm{v}_{V,i} = a(v,v)_i$.
Let $u$ be the solution to 
\eqref{eqn:primal}, 
$u_{ms}\in V_{\text{off}}$ be the solution to \eqref{eq:globalG},
$z$ the solution to \eqref{eqn:dual}, and  
$z_{ms}\in V_{\text{off}}$ the solution to \eqref{eqn:discrete_dual}.
Using Galerkin orthogonality and the relation between the primal and dual
problems, the error in the quantity of interest satisfies
\begin{align}\label{eqn:error_rep1}
g(u - u_{ms}) & = a(z,u - u_{ms}) = a(u - u_{ms}, z) = a(u - u_{ms}, z-z_{ms}).
\end{align}
Decomposing the global integration into neighborhoods by the partition of 
unity functions $\chi_i$ given by~\eqref{pou}
\begin{align}\label{eqn:error_rep2}
a(u-u_{ms},z - z_{ms}) & = \sum_{i = 1}^n \int_{\Omega} \chi_i\kappa \nabla (u-u_{ms}) \cdot \nabla (z - z_{ms}) \; dx
\le \sum_{i = 1}^n \norm{u-u_{ms}}_{V,i}\norm{z - z_{ms}}_{V,i}
\end{align}
where we used the fact that $\chi_i$ is supported in $\omega_i$
and the fact that $| \chi_i | \leq 1$.
Instead of using the norm of local residual for the primal problem defined in 
\eqref{eq:Rnorm} to be our indicator, \eqref{eqn:error_rep2} suggests
using the product of norms of 
local residuals for the primal and dual problems, posed in the same
discrete space, $V_{\text{off}}$. 
The local dual residual $R^z_i:V_i\rightarrow \mathbb{R}$ is defined by
\begin{equation}\label{eqn:dualresi}
R^z_i(v) =  g(v) - \int_{\omega_i} \kappa(x) \nabla z_{\text{ms}}\cdot \nabla v,
\end{equation}
where $z_\text{ms}\in V_{\text{off}}$ 
is the solution to  
\eqref{eqn:discrete_dual}. Analogous to \eqref{eq:Rnorm}, 
the $H^{-1}$ norm of $R^z_i$ is defined as
\begin{equation}
\| R^z_i \|_{V_i^*} = \sup_{v\in V_i} \frac{ |R^z_i(v)| }{\| v\|_{V_i}}.
\label{eq:Rnorm_z}
\end{equation}
The local version of ~\eqref{eq:a-post-bound}, applied to both primal
and dual residuals, namely
\begin{align}\label{eqn:localresbound}
\norm{u - u_{ms}}_{V,i} \le C \norm{R_i^u}_{V_i^\ast}
(\lambda_{l_i + 1}^{\omega_i})^{-1/2}, \quad
\norm{z - z_{ms}}_{V,i} \le C \norm{R_i^z}_{V_i^\ast}
(\lambda_{l_i + 1}^{\omega_i})^{-1/2},
\end{align}
motivates the local error indicator, $\eta_i$,  defined as 
\begin{align}\label{eqn:resindic}
\eta^2_i
=\| R^z_i\|_{V_i^*}  \| R^u_i\|_{V_i^*} (\lambda^{\omega_i}_{l_i+1})^{-1},
\end{align}
where $\| R^u_i\|_{V_i^*}$ and $\| R^z_i\|_{V_i^*}$ are defined in 
\eqref{eq:Rnorm} and \eqref{eq:Rnorm_z}, respectively.
Applying \eqref{eqn:resindic} and \eqref{eqn:localresbound} to 
\eqref{eqn:error_rep2} bounds the error in the goal function by
\begin{align}\label{eqn:globalestbound}
g(u - u_{ms}) \le C \sum_{i = 1}^n  \eta_i^2.
\end{align} 
In summary, the goal-error over the global domain $D$ is bounded by the 
product of energy errors of the primal and dual problems, which is in 
turn bounded by the sum of the indicators given by~\eqref{eqn:resindic},
modified by the partition of unity functions. This upper bounds suggests
the adequacy of the indicators in reducing the error.  The efficiency and
a formal convergence analysis are however not addressed here. 
As in~\cite{BuNa15} where a similar indicator is used for $hp$-refinement,
this indicator displays similar behavior to the DWR-type indicator, as shown 
in the numerical experiments; however, it is more amenable to analysis.
This indicator has the added advantage of reduced computational cost as
compared to the DWR-type indicator described below, as both primal and
dual problems are solved over the same discrete spaces, whereas for the
DWR-type method, the dual problem has greater computational complexity
than the primal.

\subsubsection*{DWR-type goal-oriented indicator}
The next error indicator is similar to DWR error indicator. 
For the primal problem solved in discrete space $V_{\text{off}}$, the 
DWR indicator is motived by the following residual equation. For 
$z$ the solution to \eqref{eqn:dual}, $u$ the solution to 
\eqref{eqn:primal}, and 
$u_{ms}\in V_{\text{off}}$ the solution to \eqref{eq:globalG} 
\begin{align}\label{eqn:dwr_sat}
g(u - u_{ms}) = a(u - u_{ms},z - z_{\text{off}})= R^u(z-z_{\text{off}}),
\end{align}
where $z_{\text{off}}$ in $V_{\text{off}}$ is arbitrary
and the global residual 
$R^u(v) = \int_D fv - \int_D \kappa(x) \nabla u_{\text{ms}}\cdot \nabla v$. 
We let $z_{\text{off}}^i$ be the component of $z_{\text{off}}$ spanned by the basis functions
corresponding to the coarse neighborhood $\omega_i$.
Localizing~\eqref{eqn:dwr_sat} by the partition of unity functions $\chi_i$,
\begin{align}\label{eqn:dwr_sat2}
R^u(z-z_{\text{off}})  = \sum_{i = 1}^N R_i^u \Big( \chi_i z-z_{\text{off}}^i \Big)
= \sum_{i = 1}^N R_i^u \Big( z_{\text{enrich}}^i-z_{\text{off}}^i \Big) + 
 \sum_{i = 1}^N R_i^u \Big( \chi_i z-z_{\text{enrich}}^i \Big). 
\end{align}
As the exact solution $z$ is unavailable, one generally instead replaces $z$
by $z_{\text{enrich}}$, 
a discrete solution from a more enriched space than the primal,
essentially neglecting the last term of \eqref{eqn:dwr_sat2}.
The function $z_{\text{enrich}}^i$ is the component of $z_{\text{enrich}}$ spanned by the basis functions
corresponding to the coarse neighborhood $\omega_i$.
In standard finite element methods the global residual is then solved
elementwise and used as an indicator.  By Galerkin orthogonality, then function
$z_{\text{off}}$ may be taken as any function in $V_{\text{off}}$ but in 
practice is taken as the projection of the enriched dual solution into 
$V_{\text{off}}$.

In this case, the dual problem is solved in the  enriched space, called
$V_\text{enrich}$.
The space $V_\text{enrich}$ is obtained by adding more basis functions 
to each coarse neighborhood. 
Recalling the construction of basis functions $\psi_{i,k}$ in 
\eqref{cgbasis},
the enriched space is constructed with more than
 $l_i$ basis functions per coarse neighborhood, specifically
\begin{equation} \label{cgbasis1}
\psi_{i,k} = \chi_i \psi_k^{\omega_i, \text{off}}, \quad \text{for} \, \, \,
1 \leq i \leq N \, \, \,  \text{and} \, \, \, 1 \leq k \leq l_i+m,
\end{equation}
where $m$ basis functions are added for each $\omega_i$.
The span of these basis is our $V_\text{enrich}$.
Let $z_\text{enrich} \in V_\text{enrich}$ be the solution for the dual problem in 
$V_\text{enrich}$, that is, $z_\text{enrich}$ satisfies 
\begin{equation}
a(v,z_{\text{enrich}}) = g(v), \;\forall v\in V_{\text{enrich}}.
\label{eqn:discrete_enrich}
\end{equation}
The DWR-type error estimator is defined as
\begin{align}\label{eqn:eta_dwr}
\eta^2_i= \Big |R^u_i \Big(P_{i}(z_\text{enrich})-\pi\big(P_{i}(z_\text{enrich})\big) \Big) \Big|, \quad \quad i=1,2,\cdots,N.
\end{align}
In the above definition, $P_i(z_{\text{enrich}})$ is the component of 
$z_{\text{enrich}}$
spanned by the basis functions $\psi_{i,k}$, $k=1,2,\cdots, l_i+m$, 
corresponding to the coarse neighborhood $\omega_i$. 
Moreover $\pi\big(P_{i}(z_\text{enrich})\big)$ is the component of 
$P_i(z_{\text{enrich}})$
spanned by the basis functions $\psi_{i,k}$, $k=1,2,\cdots, l_i$, 
in the offline space. 
Comparison with~\eqref{eqn:dwr_sat2} yields an heuristic bound, modulo the
error term created by replacing $z$ by $z_{\text{enrich}}$.


Each of the two indicators, given by \eqref{eqn:resindic}, and respectively,
\eqref{eqn:eta_dwr}, may be implemented in an adaptive framework
to determine which coarse neighborhoods to enrich.
The goal-oriented variant of the adaptive enrichment algorithm in 
Section~\ref{subsec:adaptive_algorithm} is now described.

\subsubsection*{Goal-oriented adaptive enrichment algorithm}
Choose a fixed marking parameter $0  < \theta < 1$. 
Choose also an initial offline space $V_{\text{off}}^1$ 
by specifying a fixed number of basis functions for each coarse neighborhood,
and this number is denoted by $l_i^1$.
Then, generate a sequence of spaces $V_{\text{off}}^m$
and a sequence of multiscale solutions $u_{\text{ms}}^m$ obtained by solving 
\eqref{eq:globalG}.
Specifically, for each $m=1,2,\cdots$, perform the following calculations:

\begin{enumerate}
\item[Step 1:] Find the multiscale solution in the current space. That is,
find $u_{\text{ms}}^m \in V^m_{\text{off}}$ such that
\begin{equation}
a(u^m_{\text{ms}}, v) = (f, v) \quad \text{for all} \,\,\, v \in V^m_{\text{off}}.
\label{eq:solve_u}
\end{equation}

\item[Step 2:] Find the multiscale dual solution in the current space or an enriched space. That is,
find $z_{\text{ms}}^m \in V^m_{\text{dual}}$ such that
\begin{equation}
a(z^m_{\text{ms}}, v) = (f, v) \quad \text{for all} \,\,\, v \in V^m_{\text{dual}}
\label{eq:solve_z}
\end{equation}
where 
\begin{equation*}
V^m_{\text{dual}} =
\begin{cases}
V^m_\text{off}, & \text{ for $H^{-1}$-based error estimator},\\
V^m_\text{enrich},  & \text{ for DWR-type error estimator}.
\end{cases}
\end{equation*}

\item[Step 3:] Compute the local residual. For each coarse region $\omega_i$, 
 compute
\begin{equation*}
\eta^2_i =
\begin{cases}
\| R^z_i\|_{V_i^*},  \| R^u_i\|_{V_i^*} (\lambda^{\omega_i}_{l_i+1})^{-1}, 
& \text{ for $H^{-1}$-based error estimator},\\
\big |R^u_i(P_{i}(z_\text{enrich})-\pi(P_{i}(z_\text{enrich}))) \big |,  
& \text{ for DWR-type error estimator}.
\end{cases}
\end{equation*}
\begin{align*}
R^u_i(v) &=  \int_{\omega_i} fv - \int_{\omega_i} \kappa(x)\nabla u^m_{\text{ms}}\cdot \nabla v, \\
R^z_i(v) &=  \int_{\omega_i}gv - \int_{\omega_i} \kappa(x)\nabla z^m_{\text{ms}}\cdot \nabla v,
\end{align*}
consistent with \eqref{eqn:primalresi} and~\eqref{eqn:dualresi};
and the norm is defined in (\ref{eq:Rnorm}) and (\ref{eq:Rnorm_z}) respectively.
Next, re-enumerate the coarse neighborhoods so that the above local residuals $\eta_i^2$ are arranged in decreasing order
$\eta^2_1 \geq \eta^2_2 \geq \cdots \geq \eta^2_N$.
That is, in the new enumeration,
the coarse neighborhood $\omega_1$ has the largest residual $\eta_1^2$.
As in Remark~\ref{remark:bin} the full sort of the estimators can be
replaced by a partial sort for a marked set of quasi-optimal cardinality.

\item[Step 4:] Find the coarse regions where enrichment is needed. Choose 
the smallest integer $k$ such that
\begin{equation}
\theta \sum_{i=1}^N \eta_i^2 \leq \sum_{i=1}^k \eta_i^2.
\label{eq:criteria}
\end{equation}
The coarse neighborhoods $\omega_1,\omega_2,\cdots, \omega_k,$ are then
enriched with additional basis functions. Alternately, a set of neighborhoods
based on the binning strategy for a partial
sort can be chosen as described in Step 3 of the Adaptive 
algorithm \ref{subsec:adaptive_algorithm}. 

\item[Step 5:] Enrich the space. 
For each $i=1,2,\cdots, k$, add basis functions
for the region $\omega_i$ according to the following rule.
Let $s$ be the smallest positive integer such that
$\lambda_{l_i^m+s+1}$ is large enough compared with $\lambda_{l_i^m+1}$.
Then include the eigenfunctions 
$\Psi^{\text{off}}_{l_i^m+1}, \cdots, \Psi^{\text{off}}_{l_i^m+s}$
in the construction of the basis functions.
The resulting space is denoted as $V_{\text{off}}^{m+1}$.
Mathematically, the space $V_{\text{off}}^{m+1}$ is defined as
\begin{equation*}
V_{\text{off}}^{m+1} = V_{\text{off}}^{m} + \text{span} \cup_{i=1}^k \cup_{j=l_i^m+1}^{l_i^m+s} \{  \psi_{i,j} \}
\end{equation*}
where $\psi_{i,j} = \chi_i \psi_j^{\omega_i,\text{off}}$
and $\psi_j^{\omega_i,\text{off}} = \sum_{r=1}^{l_i} \Psi_{jr}^{\text{off}} \psi_r^{\text{snap}}$,
with $j=l_i^m+1, \cdots, l_i^m+s$,
denote the new basis functions obtained by the eigenfunctions $\Psi^{\text{off}}_{l_i^m+1}, \cdots, \Psi^{\text{off}}_{l_i^m+s}$.
In addition, set $l_i^{m+1} = l_i^m+s$.
\end{enumerate}

In the next section, we demonstrate the efficiency of the goal-oriented 
adaptive algorithm defined above on a problem with high-contrast
multiscale coefficients.  
The results are compared with the standard residual-based adaptive
method defined in the previous section.  We note both the increased 
efficiency in the reduction in goal-error, $|g(u - u_{ms})|$, and 
the decreased reduction in the energy-norm error with the goal-oriented methods.
These two observations suggest the method does what it was designed to do: 
focus the adaptive enrichment towards reduction in goal-error without
resolving the solution where it has limited influence on the goal-error.
We also note similarity in performance between the two indicators, both
demonstrating errors with a similar observed rate of convergence.
\section{Numerical Results}
\label{sec:numresults}
\begin{figure}[!th]
\centering
\includegraphics[scale=0.5]{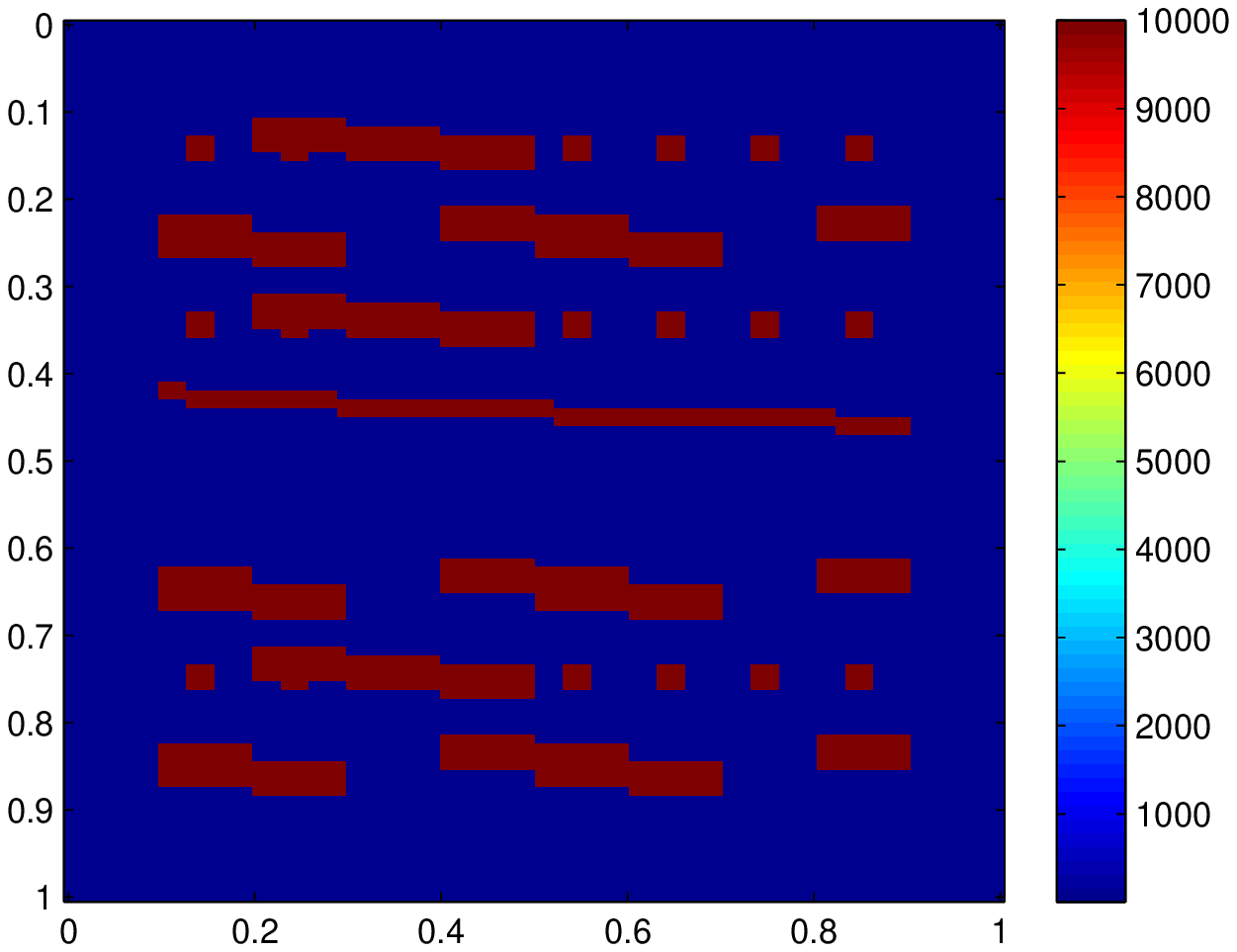}~\hfil~
\includegraphics[scale=0.5]{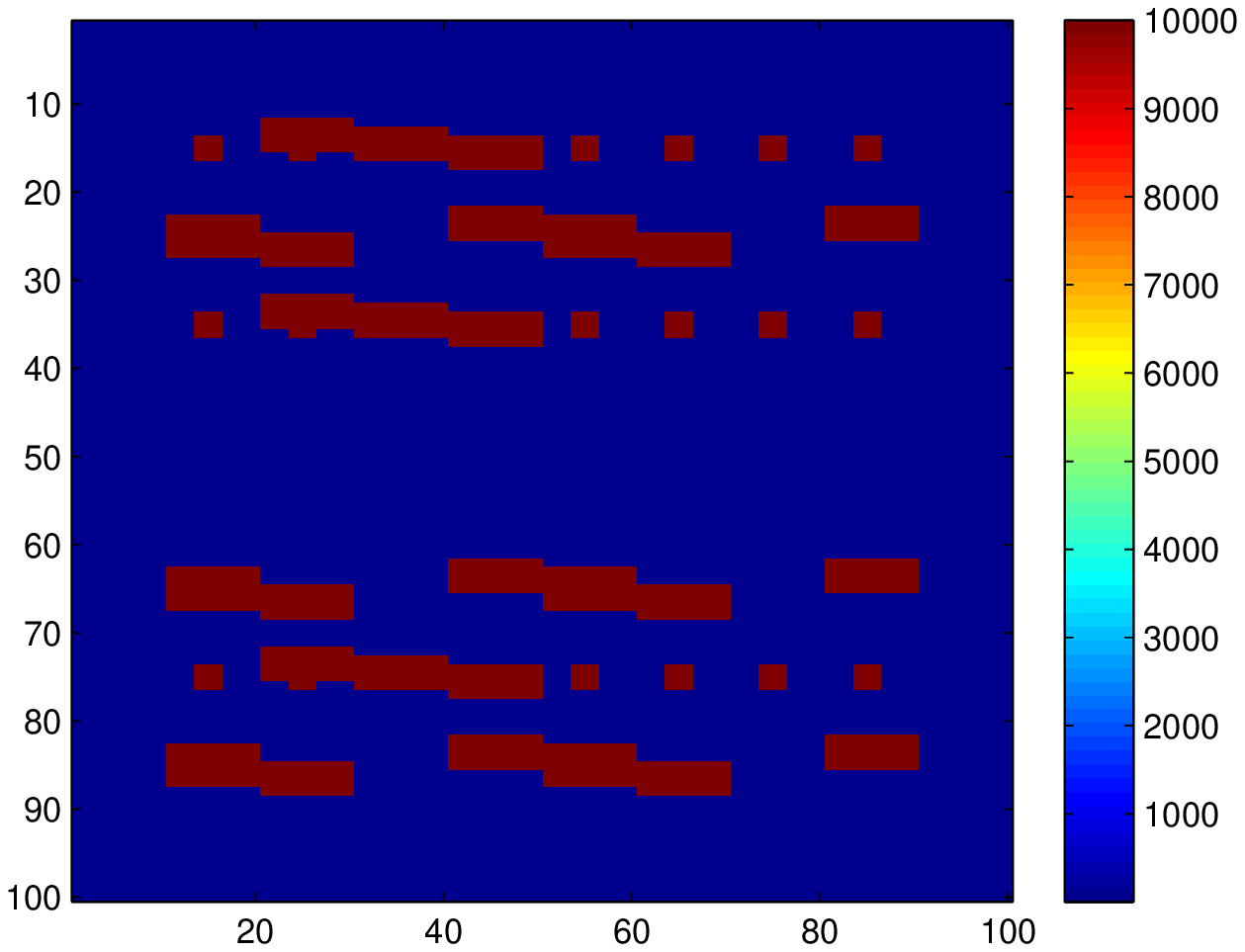}
\protect\caption{Left: the  coefficient $\kappa_1$, corresponding to
Figure \ref{fig:Error_case1}.
Right: the coefficient $\kappa_2$, corresponding to Figure 
\ref{fig:Error_case3}.}
\label{fig:medium1}
\end{figure}

%
In this section, we present two numerical examples for multiscale problems
with high-contrast coefficients and
compare the performance of the two indicators defined in the previous section.
For our simulations, we take the domain $\Omega=(0,1)^{2}$, and
the inflow-outflow source term $f=\chi_{K_{1}}-\chi_{K_{2}}$,
where $K_{1}=[0.1,0.2]\times[0.8,0.9]$, 
and $K_{2}=[0.8,0.9]\times[0.1,0.2]$.
We consider the problem of finding $g(u)$ for $u$ the solution to 
\eqref{eq:original}, namely
\begin{align}\label{eqn:exampleprob}
-\div(\kappa(x) \nabla u) = f, ~\text{ in }~D, \quad u = 0 ~\text{ on }~ \pa D.
\end{align} 
The goal functional 
\[
g(u)=\int_{K_{2}}u,
\]
is the average value of $u$ on the outflow region $K_2$.
In practice, $K_2$ is the location of the wells, and it is important that 
the average pressure $u$ on $K_2$ is accurate. 
The two examples differ by the high-contrast coefficients $\kappa(x)$, shown
in Figure~\ref{fig:medium1}. Shown on the left,
$\kappa_1$ features a high-conductivity
channel crossing the domain separating the inflow and outflow; on 
the right, $\kappa_2$ is a similar coefficient without the channel. 
These coefficents are visualized with the blue region indicating the value $1$
and the red region indicating the contrasts. 
In each example, we consider two different contrast strengths, 
$10^4$ and $10^6$. 
We note the invariance in the relative performance of the indicators 
with respect to the contrast strengths.

\begin{figure}[!th]
\centering
\includegraphics[scale=0.4]{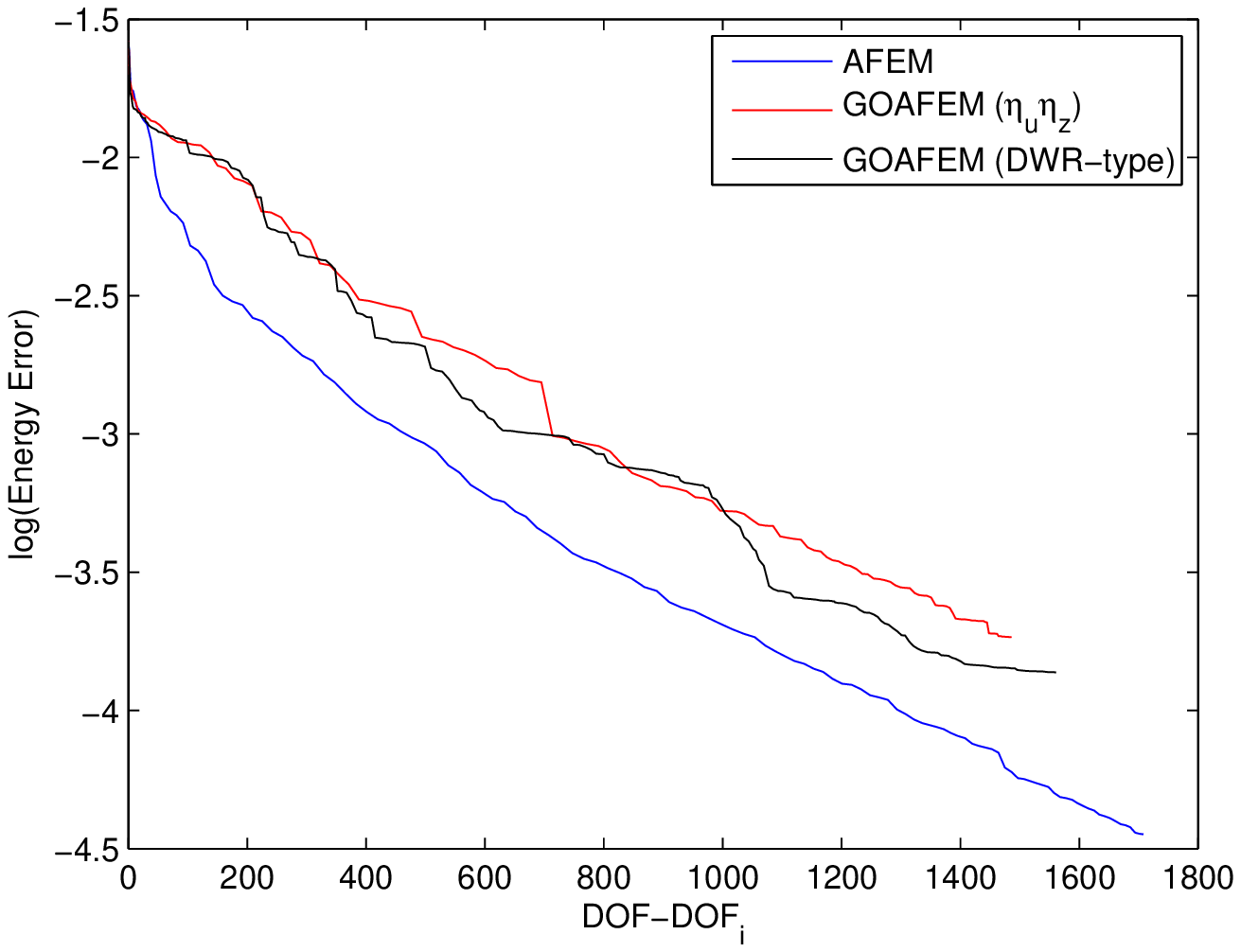} 
\includegraphics[scale=0.4]{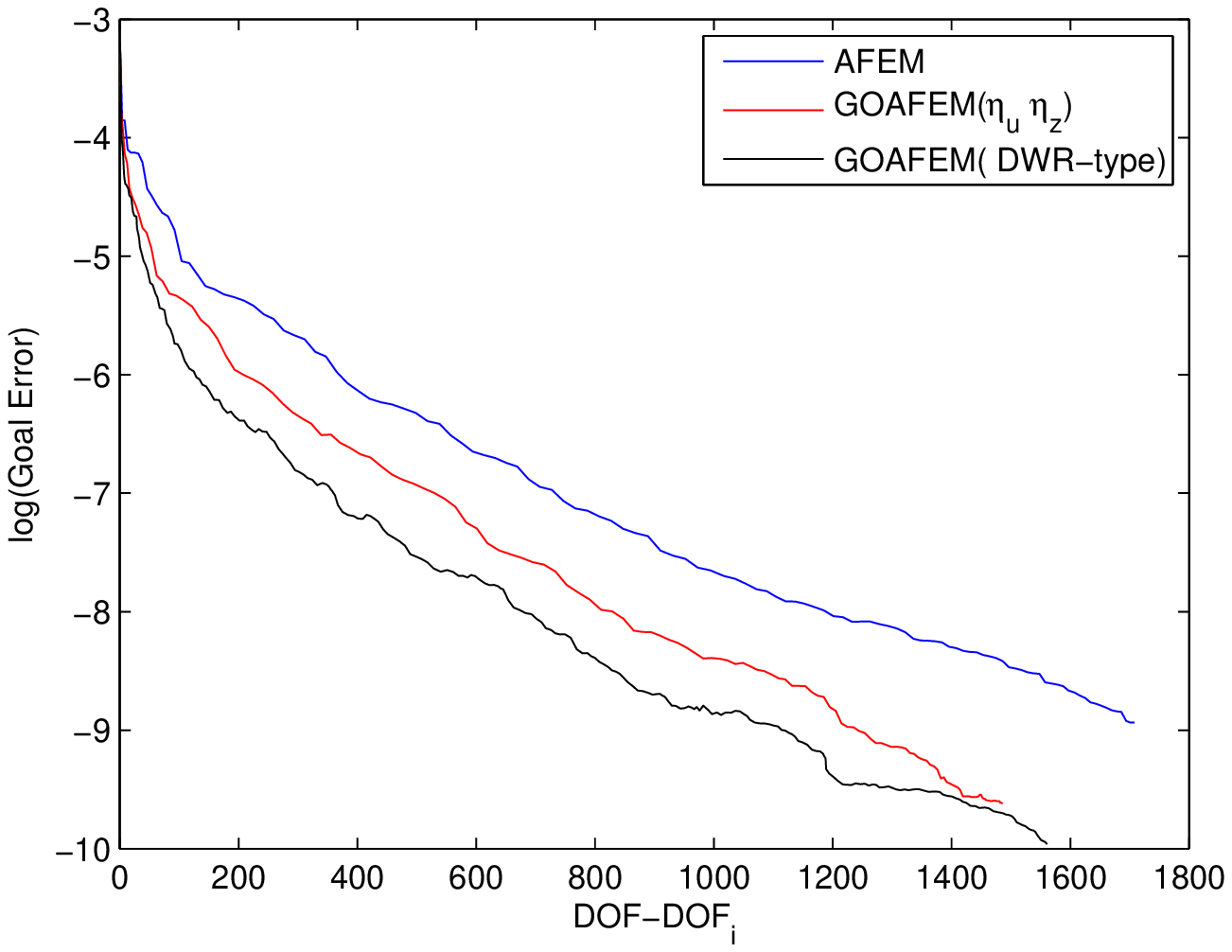}
\includegraphics[scale=0.4]{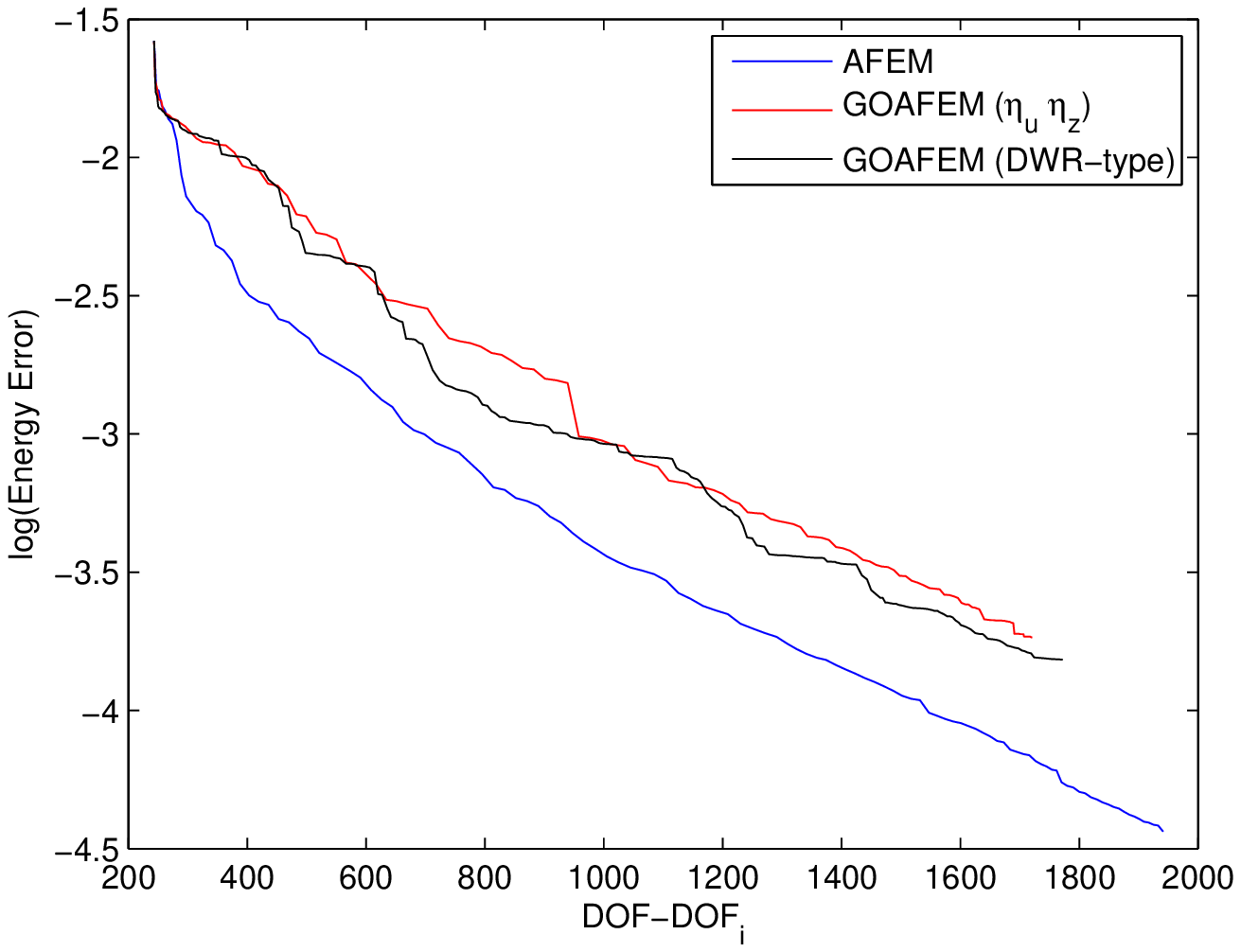} 
\includegraphics[scale=0.4]{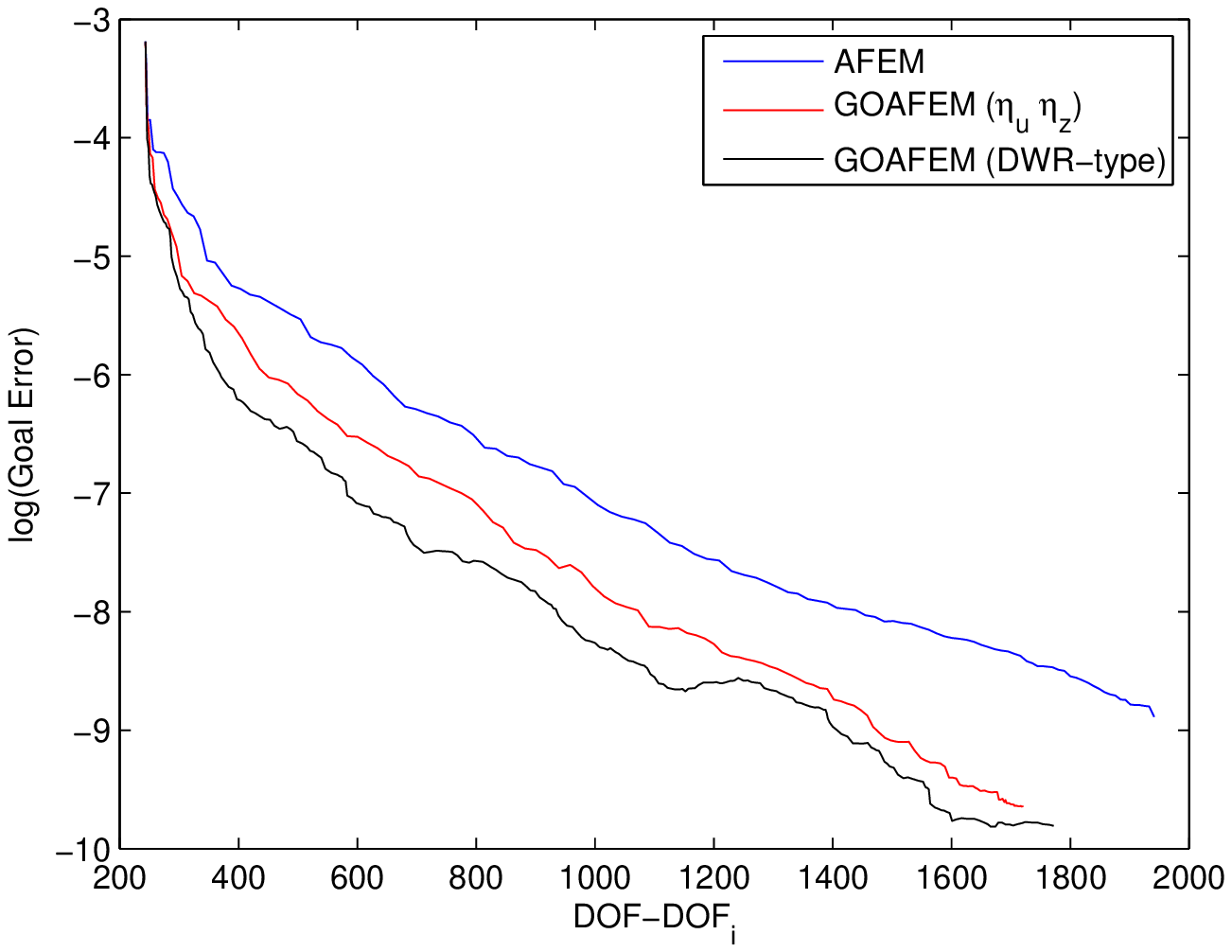}
\protect\caption{Lower contrast ($10^4$). 
Top left:  $\log \|u-u_{ms}\|_V$. 
Top right: $\log |g(u-u_{ms})|$. 
\protect \\ Higher contrast ($10^6$). 
Bottom left:  $\log \|u-u_{ms}\|_V$, 
Bottom right: $\log |g(u-u_{ms})|$.}
\label{fig:Error_case1}
\end{figure}

\begin{figure}[!th]
\centering
\includegraphics[scale=0.4]{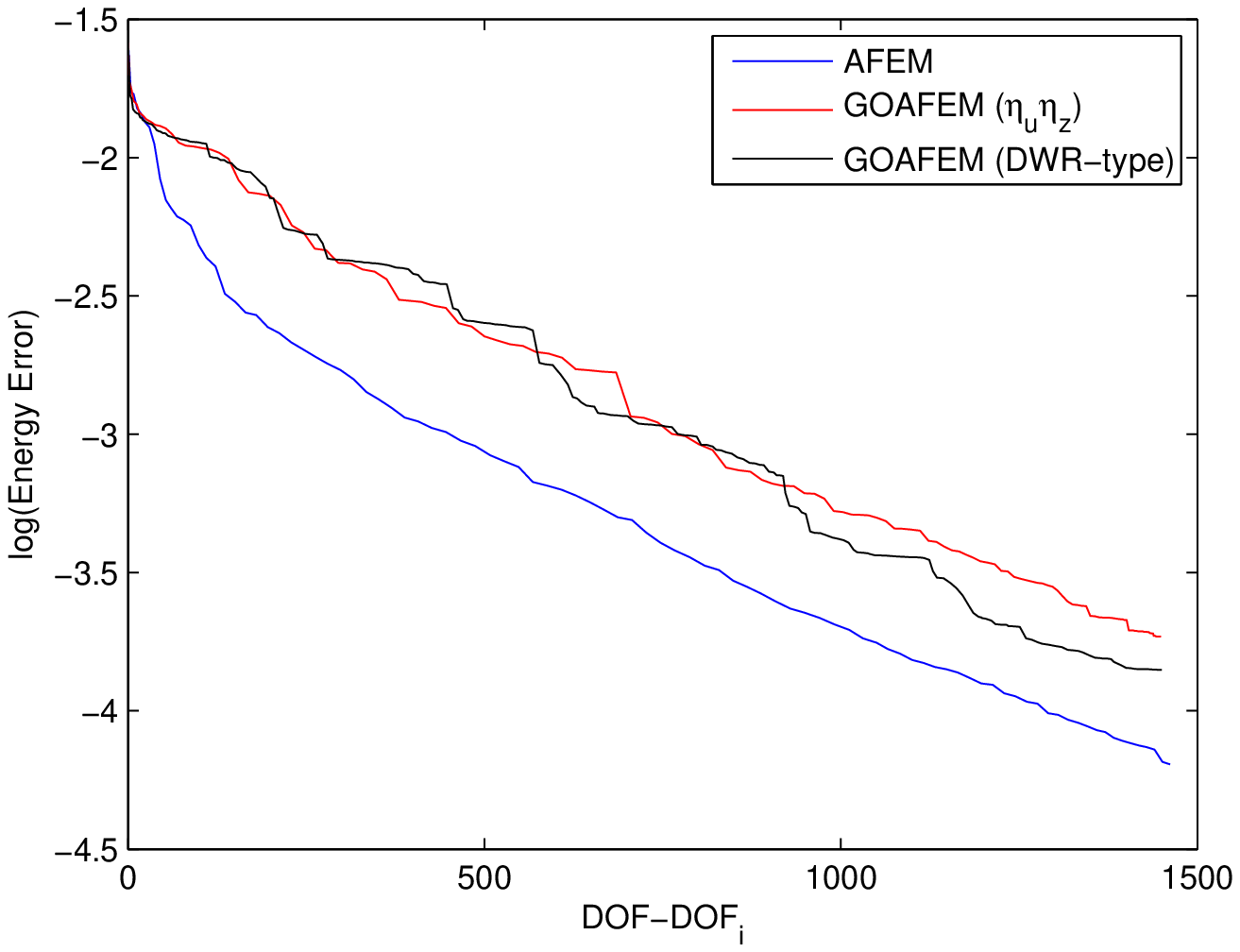} 
\includegraphics[scale=0.4]{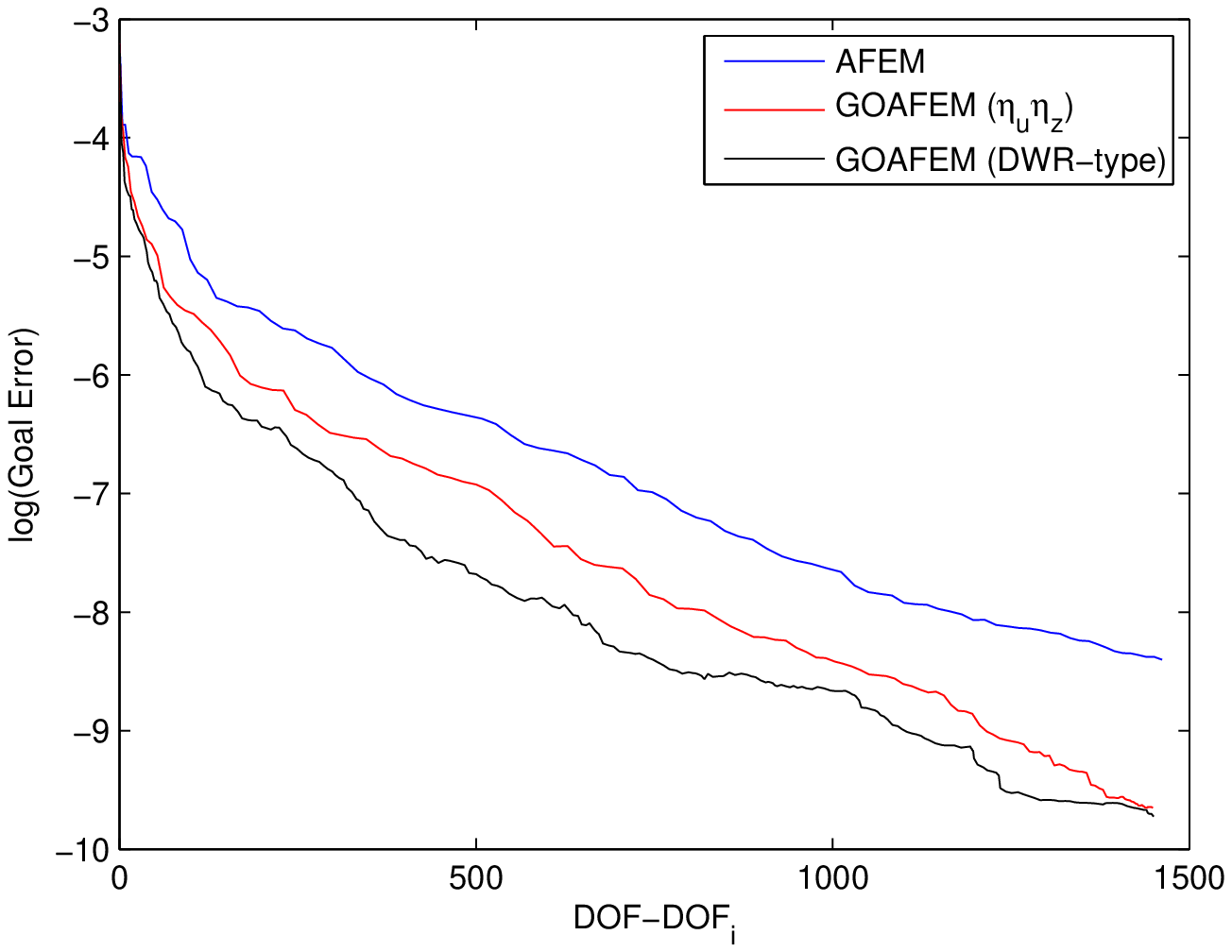}
\includegraphics[scale=0.4]{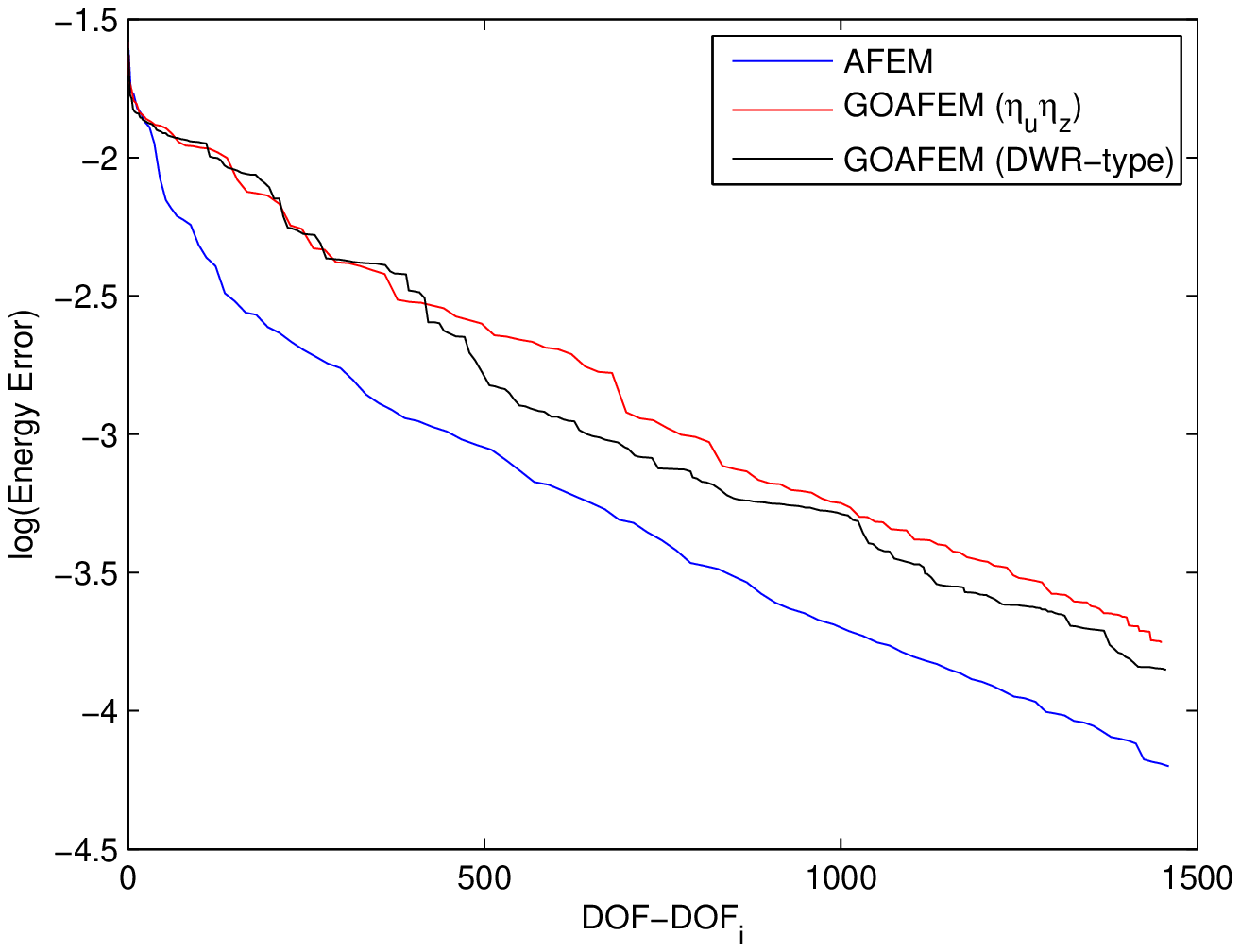} 
\includegraphics[scale=0.4]{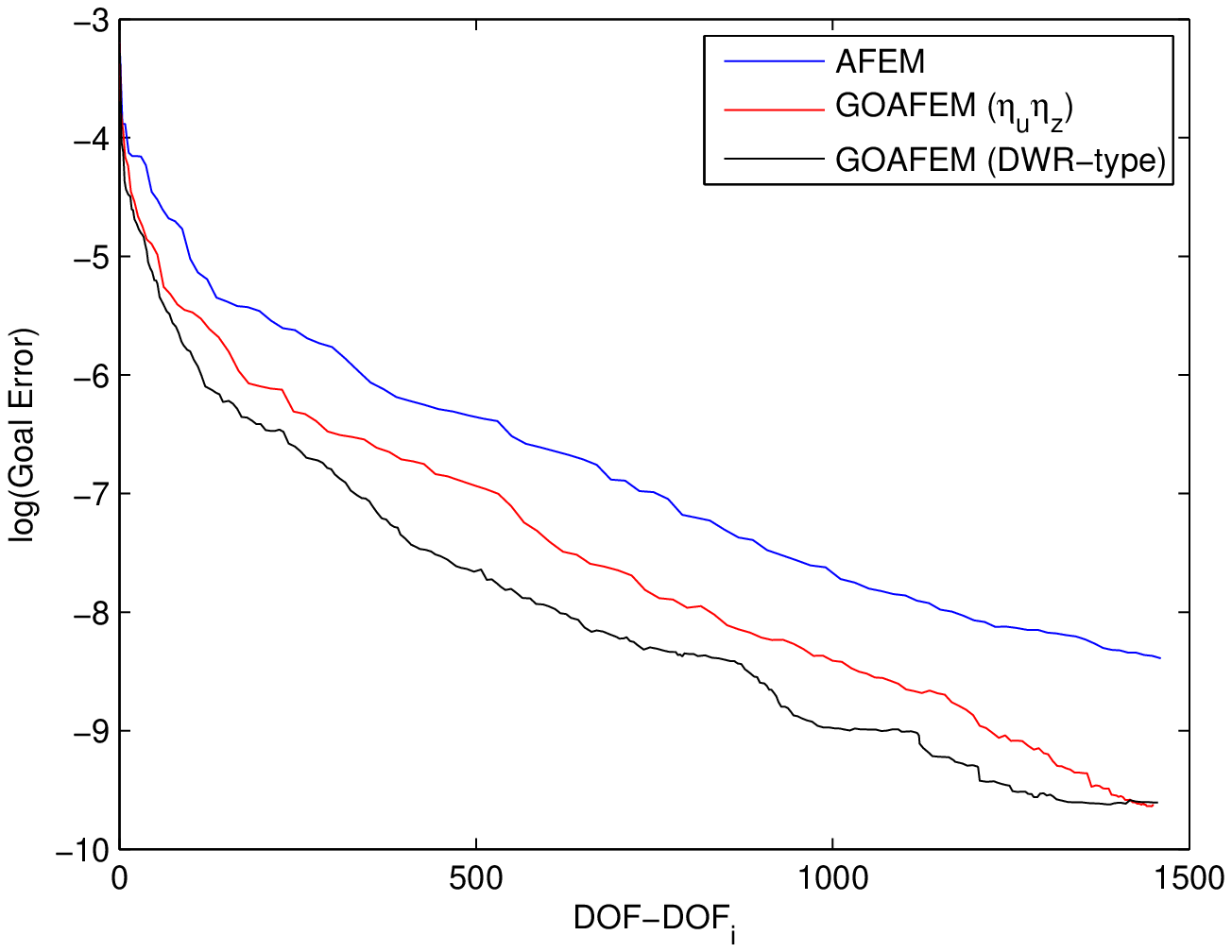}
\protect\caption{Lower contrast ($10^4$). 
Top left:  $\log \|u-u_{ms}\|_V$, 
Top right: $\log |g(u-u_{ms})|$. 
\protect \\ Higher contrast($10^6$), 
Bottom left:  $\log \|u-u_{ms}\|_V$, 
Bottom right: $\log |g(u-u_{ms})|$.}
\label{fig:Error_case3}
\end{figure}

For the first example, 
Figure \ref{fig:Error_case1} shows a comparison in error reduction between
the three indicators for coefficient $\kappa_1$. 
In the two figures on the right, 
we compare the logarithms of the energy norm errors
$g(u-u_{\text{ms}})$ against the number of unknowns 
for three types of adaptive enrichment algorithms;
namely, the residual based method described in 
Section \ref{subsec:adaptive_algorithm} 
(denoted in blue in Figure \ref{fig:Error_case1}),
the $H^{-1}$ residual-based goal-oriented method 
(denoted in red in Figure \ref{fig:Error_case1}), 
and the DWR type goal-oriented method 
(denoted in black in Figure \ref{fig:Error_case1})
as described in Section \ref{sec:errorindicator}.
From these results, we see the two types of goal-oriented methods 
behave similarly and outperform the standard adaptive method with an
improved rate of goal-error reduction. 
We note a more stable decrease in error reduction for
the goal-oriented 
residual-type method, but slightly improved, if less predictable error 
reduction for the DWR-type.  This last observation 
is to be expected, as the DWR-type 
indicator does not account for the error created by using the enriched solution
$z_{\text{enrich}}$ in place of the exact dual solution $z$, as in 
\eqref{eqn:dwr_sat2}.

On the left of Figure~\ref{fig:Error_case1}, 
we see the standard residual based method outperforms the goal-oriented 
methods
for the energy norm error, $\| u - u_{\text{ms}}\|_V$.  
This confirms that the goal-oriented methods
are driving the adaptivity toward a more efficient evaluation of the 
quantity of interest without expending additional computational effort
resolving features of the solution with limited influence on the goal.

For the second example, we consider finding $g(u)$ for $u$ that
satisfies ~\eqref{eqn:exampleprob}, with $\kappa(x)$ given by $\kappa_2$ 
shown on the 
right of Figure~\ref{fig:medium1}, with no high-conductivity channel.
As seen in Figure~\ref{fig:Error_case3}, the results are qualitatively
similar to the results of the first example.  In summary, the plots on 
the right show the goal-error reduction for the 
higher and respectively lower contrast cases. The residual-type and DWR-type
goal-oriented indicators achieve a better rate of error reduction than
the standard adaptive method, with the DWR-type showing generally the
lowest error, with the least-steady decrease. The plots on the right
of Figure ~\ref{fig:Error_case3} show the reduction of energy error for the
three indicators.  As in the first example, the standard $H^{-1}$ residual
based adaptive method designed to reduce the energy error shows the best
performance here, whereas the goal-oriented methods yield steady error 
reduction but are focused on localized error reduction in the region of
the goal-functional, rather than across the entire domain.

These results demonstrate the importance of goal-oriented adaptivity in
GMsFEM, particularly in cases where the global domain is significantly
larger than the region of infuence for the quantity of interest.  In
partiuclar, problems with many localized features only some of which 
significantly influence the quantity of interest will benefit 
from a goal-oriented adaptive strategy. 

\section{Conclusion}
In this paper we defined two types of error indicators that can be used
in an adaptive algorithm for multiscale problems with high-contrast 
coefficients.  The goal-oriented adaptive algorithm fits within the
framework of GMsFEM, and focuses the adaptivity on reducing the error
in the quantity of interest, rather than in global norm. 
We first reviewed the general ideas of GMsFEM for high-contrast problems, 
then gave a detailed overview of the construction of multiscale basis 
functions, and a residual based adaptive algorithm designed to reduce
the energy norm error.  We stated the dual problem, 
then motivated and introduced two goal-oriented 
error indicators, and described their use in a a goal-oriented adaptive 
algorithm.  Finally, we demonstrated the efficiency of the goal-oriented
algorithm and estimators compared with the standard adaptive method
introduced earlier.  We found for both indicators, the goal-oriented 
method reduced the error in the goal-function at a better rate than the
standard method.  We also found the two indicators perform similarly, with 
some increase in error reduction seen in the DWR-type indicator, but at
the cost of solving the dual problem in a more enriched space, increasing
the computational complexity. The residual-based indicator on the other
hand may be more amenable to convergence analysis, as may be investigated
in future work. The current results indicate the goal-oriented strategy
increases the efficiency of GMsFEM when a function of the solution rather
than the solution in its entirety is of interest.

\bibliographystyle{plain}
\bibliography{references1}
\end{document}